\documentclass{article}
\usepackage{graphicx,amsmath,amssymb,amsthm}

\numberwithin{equation}{section}

\theoremstyle{plane}
\newtheorem{thm}{Theorem}[section]

\newtheorem{prop}[thm]{Proposition}
\newtheorem{cor}[thm]{Corollary}

\theoremstyle{remark}
\newtheorem*{rem}{Remark}

\newcommand{\vol}{\mathop{\mathrm{Vol}}\nolimits}
\newcommand{\li}{\mathop{\mathrm{Li_2}}\nolimits}

\title{A volume formula for hyperbolic tetrahedra in terms of edge lengths}
\author{Jun {\sc Murakami} \vspace{\medskipamount}\\
Department of Mathematical Science\\
School of Science and Engineering\\
Waseda University\\
Tokyo 169-8555, JAPAN\\
{\tt murakami@waseda.jp}
\vspace{\bigskipamount}
\and Akira {\sc Ushijima} \vspace{\medskipamount}\\
Mathematics Institute\\
University of Warwick\\
Coventry CV4 7AL, UK\\
{\tt ushijima@maths.warwick.ac.uk}
}

\date{\empty}

\begin{document} \maketitle
\begin{abstract}
We give a closed formula for volumes of generic hyperbolic tetrahedra
in terms of edge lengths.
The cue of our formula is  by the
volume conjecture for the Turaev-Viro invariant of closed
$3$-manifolds, which is defined from the quantum $6j$-symbols.
This formula contains the dilogarithm functions,
and we specify the adequate branch to get the actual value
of the volumes.
\end{abstract}
\begin{flushleft}\begin{small}
{\bf Keywords:}
hyperbolic tetrahedron, quantum $6j$-symbol, volume formula.\\
{\bf 2000 Mathematics Subject Classifications:}
primary: 52A38; secondary: 51M09.
\end{small}\end{flushleft}

\section{Introduction}
To determine the volume of a hyperbolic
tetrahedron is one of the fundamental
problem of the hyperbolic $3$-space.
For orthoschemes, a closed
formula is known at the very
beginning of the hyperbolic geometry.
However, closed formulas for the volumes of
generic hyperbolic tetrahedra are given rather recently by
\cite{ck} and is modified to more symmetric form in \cite{my}.
These formulas are given in terms of dihedral angles of the tetrahedron.

In this paper, we give a closed formula in terms of edge lengths 
instead of dihedral angles.
The background of this formula is the volume conjecture for hyperbolic $3$-manifolds.
In \cite{ka}, R. Kashaev conjectured in a relation 
between certain quantum invariants of hyperbolic knots 
and the volumes of the knot complements of them, 
and checked it numerically for some examples.
In \cite{mu}, H. Murakami observed that such relation also holds 
for some closed hyperbolic $3$-manifolds.
These observations suggest that there is a relation 
between the volume of a hyperbolic $3$-manifold and its Turaev-Viro invariant (see \cite{tv}),
which is a state-sum invariant given by quantum $6j$-symbols 
associated to the tetrahedra of a simplicial decomposition of the manifold.
Indeed, applying Kashaev's method to investigating the asymptotics of the quantum $6j$-symbols, 
we get the formula for the volume of a hyperbolic tetrahedron.
The formula by dihedral angles in \cite{my} also comes from the $6j$-symbols in the same way, 
and it resembles so much the formula given in this paper.

Such similarity comes from
the relations in \cite[p.\ 294]{mi} and \cite[p.\ 311]{sa}, by which
the volume function
in terms of edge lengths
implicitly comes from
the volume function
in terms of dihedral angles.
We can apply this method
to the formula in \cite{my}.
However, to get the actual
volumes, we have to specify
the adequate branches of the
dilogarithm functions in
the formula.
Here we give a direct proof of
our formula,
 and show that we can
actually compute
 the volume of
any hyperbolic tetrahedron
from it.
\bigskip\par\noindent
{\bf Acknowledgement.}
We would like to thank Professor Yana Mohanty for informing the
relations in \cite{mi} and \cite{sa}.
The second author would like to thank the support 
from the JSPS Postdoctoral Fellowships for Research Abroad $2003$.

%
%
%
%
%Section 2
%
%
%
%
\section{Volume formulae for hyperbolic tetrahedra}
\label{sec_formulae}

For complex numbers $z , a_1 , a_2 , \ldots , a_6$,
we first define a complex-valued function $U=U (a_1, a_2, a_3, a_4, a_5, a_6, z)$ 
as follows:
\begin{multline*}
U := \li (z) + \li (a_1 a_2 a_4 a_5 z) 
+ \li (a_1 a_3 a_4 a_6 z) + \li (a_2 a_3 a_5 a_6 z)\\
 - \li (- a_1 a_2 a_3 z) - \li (- a_1 a_5 a_6 z) 
- \li (- a_2 a_4 a_6 z) - \li (- a_3 a_4 a_5 z) \, ,
\end{multline*}
where $\li ( \cdot )$ is the dilogarithm function 
defined by the analytic continuation of the following integral:
\begin{equation*}
\li (x) := - \int_{0}^{x} \frac{\log (1-t )}{t} \mathrm{d} t 
	\quad \text{ for a real number $x < 1$.}
\end{equation*}

Secondly we define two complex values $z_-$ and $z_+$ as follows:
\begin{equation*}
z_- := \frac{-q_1 - \sqrt{q_1^2 - 4 \, q_0 \, q_2}}{2 \, q_2} \quad \text{and} \quad
z_+ := \frac{-q_1 + \sqrt{q_1^2 - 4 \, q_0 \, q_2}}{2 \, q_2} \, ,
\end{equation*}
where
\begin{align*}
q_0 &:=
		1 + {a_1}\,{a_2}\,{a_3} + {a_1}\,{a_5}\,{a_6} + 
		{a_2}\,{a_4}\,{a_6} + {a_3}\,{a_4}\,{a_5}\\
	&\qquad + {a_1}\,{a_2}\,{a_4}\,{a_5} + 
		{a_1}\,{a_3}\,{a_4}\,{a_6} + 
		{a_2}\,{a_3}\,{a_5}\,{a_6} \, ,\\
q_1 &:=
		- {a_1}\,{a_2}\,{a_3}\,{a_4}\,{a_5}\,{a_6}\,
		\left\{ \left( a_1 - a_1^{-1} \right) \left( a_4 - a_4^{-1} \right) \right.\\
	&\qquad \left. + 
		\left( a_2 - a_2^{-1} \right) \left( a_5 - a_5^{-1} \right)  + 
		\left( a_3 - a_3^{-1} \right) \left( a_6 - a_6^{-1} \right) \right\} ,\\
q_2 &:=
		{a_1}\,{a_2}\,{a_3}\,{a_4}\,{a_5}\,{a_6}\, \left( 
		{a_1}\,{a_4} + {a_2}\,{a_5} + {a_3}\,{a_6} \right.\\
	&\qquad \left. + {a_1}\,{a_2}\,{a_6} + 
		{a_1}\,{a_3}\,{a_5} + 
		{a_2}\,{a_3}\,{a_4} + 
    	{a_4}\,{a_5}\,{a_6} + 
    	{a_1}\,{a_2}\,{a_3}\,{a_4}\,{a_5}\,{a_6} \right) .
\end{align*}
These $z_-$ and $z_+$ are originally defined in \cite{my} 
as two non-trivial solutions of 
$\displaystyle{ \Re ( z \, \frac{\partial U}{\partial z} ) = 0 }$,
where $\Re$ means the real part.
We will see it (and other properties of $z_-$ and $z_+$) in Section~\ref{s_proof}.

\begin{figure}[ht]
        \begin{center}
        \includegraphics[clip]{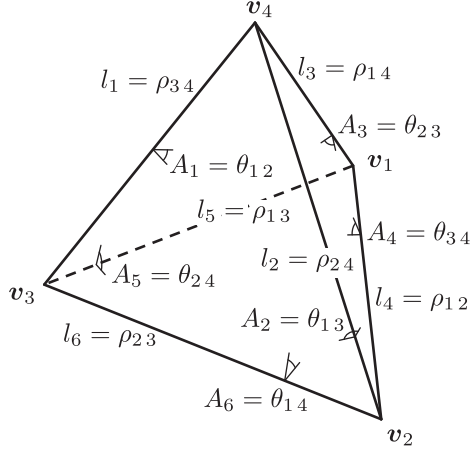}  
        \end{center}
        \caption{Dihedral angles and edge lengths of a tetrahedron}
        \label{f_tetra}
\end{figure}

Finally we define a complex-valued function
$V ( a_1, a_2, a_3, a_4, a_5, a_6 )$ as follows:
\begin{multline*}
V ( a_1, a_2, a_3, a_4, a_5, a_6 )\\
	 := \frac{\sqrt{-1}}{4} \left\{
		\left( U (a_1, a_2, a_3, a_4, a_5, a_6, z_-) - z_- 
		\left. \frac{\partial U}{\partial z} \right \rvert_{z = z_-} \log z_- \right) \right.\\
	\left. - \left( U (a_1, a_2, a_3, a_4, a_5, a_6, z_+) - z_+
		\left. \frac{\partial U}{\partial z} \right \rvert_{z = z_+} \log z_+ \right) \right\} .
\end{multline*}
Then the following theorem has been proved in \cite{my}:
%
%
%
%
%Theorem 2.1
%
%
%
%
\vspace{\medskipamount} \begin{thm} \label{thm_va}
Let\/ $T=T ( A_1, A_2, A_3, A_4, A_5, A_6 )$ be a %hyperbolic or spherical 
tetrahedron with dihedral angles as shown in Figure~\ref{f_tetra},
and let\/ $\vol (T)$ be its volume.
Then the following equations hold:
\begin{enumerate}
\item if\/ $T$ is hyperbolic, then $\vol (T) = - V_A$,
\item if\/ $T$ is Euclidean, then $V_A = 0$, and
\item if\/ $T$ is spherical, then $\vol (T) = - \sqrt{-1} \, V_A$,
\end{enumerate}
where\/ $V_A := V (e^{\sqrt{-1} \, A_1}, e^{\sqrt{-1} \, A_2}, e^{\sqrt{-1} \, A_3}, 
e^{\sqrt{-1} \, A_4}, e^{\sqrt{-1} \, A_5}, e^{\sqrt{-1} \, A_6}) $.
\end{thm}
We can easily see that both $z_-$ and $z_+$ have the following simpler presentations in this case:
\begin{align*}
z_- &= - \frac{2}{q_A} \left(
	\sin A_1 \, \sin A_4 + \sin A_2 \, \sin A_5 + \sin A_3 \, \sin A_6 - \sqrt{\det G_A} 
	\right) ,\\
z_+ &= - \frac{2}{q_A} \left(
	\sin A_1 \, \sin A_4 + \sin A_2 \, \sin A_5 + \sin A_3 \, \sin A_6 + \sqrt{\det G_A} 
	\right) ,\\
\end{align*}
where
\begin{align*}
q_A &:= 
	e^{\sqrt{-1} \left( A_1 + A_4 \right)} + e^{\sqrt{-1} \left( A_2 + A_5 \right)} 
	+ e^{\sqrt{-1} \left( A_3 + A_6 \right)}\\
	&\text{ } + e^{\sqrt{-1} \left( A_1 + A_2 + A_6 \right)} 
	+ e^{\sqrt{-1} \left( A_1 + A_3 + A_5 \right)} 
	+ e^{\sqrt{-1} \left( A_2 + A_3 + A_4 \right)} 
	+ e^{\sqrt{-1} \left( A_4 + A_5 + A_6 \right)}\\
	&\text{ } + e^{\sqrt{-1} \left( A_1 + A_2 + A_3 + A_4 + A_5 + A_6 \right)}
\end{align*}
and $G_A$ is the Gram matrix
in terms of dihedral angles defined as follows 
(precise definition will be seen in Section~\ref{s_gram}):
%
%
%eq_ga
%
%
\begin{equation}
G_A :=
\begin{pmatrix}
1 & -\cos A_1 & -\cos A_2 & -\cos A_6 \\
-\cos A_1 & 1 & -\cos A_3 & -\cos A_5 \\
-\cos A_2 & -\cos A_3 & 1 & -\cos A_4 \\
-\cos A_6 & -\cos A_5 & -\cos A_4 & 1
\end{pmatrix} . \label{eq_ga}
\end{equation}
The original proof of Theorem~\ref{thm_va}
was to show the equivalence to Cho-Kim's formula in \cite{ck},
and now there are several different proofs and formulations (see \cite{dl,us}).

\vspace{\bigskipamount}
The main purpose of this paper is to give a proof of the following theorem :
%
%
%
%
%Theorem 2.2
%
%
%
%
\begin{thm} \label{thm_vl}
Let\/ $l_i$ be the length of the edge corresponding to a dihedral angle\/ $A_i$
(see Figure~\ref{f_tetra} again), and we
denote by\/ $T=T ( l_1, l_2, l_3, l_4, l_5, l_6 )$ a hyperbolic tetrahedron 
with edge lengths\/ $l_i$.
Then the following equation holds:
\begin{equation*}
\vol (T) = V_l - \sum_{i = 1}^6 l_i \, \frac{\partial V_l}{\partial l_i},
\end{equation*}
where\/ $V_l := V (-e^{l_4}, -e^{l_5}, -e^{l_6}, -e^{l_1}, -e^{l_2}, -e^{l_3} ) $ .
\end{thm} \vspace{\medskipamount} 
We here note that, in this case, both $z_-$ and $z_+$ have the following simpler presentations:
\begin{align*}
z_- &= \frac{2}{q_l} \left(
	\sinh l_1 \, \sinh l_4 + \sinh l_2 \, \sinh l_5 + \sinh l_3 \, \sinh l_6 - \sqrt{\det G_l} 
	\right) ,\\
z_+ &= \frac{2}{q_l} \left(
	\sinh l_1 \, \sinh l_4 + \sinh l_2 \, \sinh l_5 + \sinh l_3 \, \sinh l_6 + \sqrt{\det G_l} 
	\right) ,\\
\end{align*}
where 
\begin{multline*}
q_l := e^{\left( l_1 + l_4 \right)} + e^{\left( l_2 + l_5 \right)} 
	+ e^{\left( l_3 + l_6 \right)}\\
	- e^{  \left( l_1 + l_2 + l_3 \right)}
	- e^{  \left( l_1 + l_5 + l_6 \right)}
	- e^{  \left( l_2 + l_4 + l_6 \right)}
	- e^{  \left( l_3 + l_4 + l_5 \right)}
	+ e^{  \left( l_1 + l_2 + l_3 + l_4 + l_5 + l_6 \right)}
\end{multline*}
and $G_l$ is the Gram matrix
in terms of edge lengths defined as follows
(precise definition will be seen in Section~\ref{s_gram}):
%
%
%eq_gl
%
%
\begin{equation}
G_l := \begin{pmatrix}
-1 & -\cosh l_4 & -\cosh l_5 & -\cosh l_3 \\
-\cosh l_4 & -1 & -\cosh l_6 & -\cosh l_2 \\
-\cosh l_5 & -\cosh l_6 & -1 & -\cosh l_1 \\
-\cosh l_3 & -\cosh l_2 & -\cosh l_1 & -1
\end{pmatrix} . \label{eq_gl}
\end{equation}

%
%
%
%
%Section 3
%
%
%
%
\section{A relationship between dihedral angles and edge lengths 
of a hyperbolic simplex}
\label{s_gram}
Unlike Euclidean one, in hyperbolic geometry
the dihedral angles and the lengths of the edges of a simplex is closely related each other.
In this section we review two well-known aspects of the relationship, namely
we review how to obtain dihedral angles from edge lengths, and vice verse.
See, for example, \cite{avs, lu} for details.

For any $n \geq 3$,
let $\left( \mathbb{H}^n , \left \langle \, \cdot ,  \cdot \, \right \rangle \right)$
be the {\em hyperboloid model of an $n$-dimensional hyperbolic space}.
For a hyperbolic $n$-simplex $\triangle \subset \mathbb{H}^n$,
we denote by $\left\{\boldsymbol{v}_i \right\}_{i = 1}^{n+1}$
the set of vectors of the vertices,
and by $\rho_{i \, j}$ the length of the edge 
with end points $\boldsymbol{v}_i$ and $\boldsymbol{v}_j$.
Then the following relation holds:
\begin{equation} \label{eq_iplength}
\left \langle \boldsymbol{v}_i , \boldsymbol{v}_j \right \rangle 
= - \cosh \rho_{i \, j} \, .
\end{equation}

We denote by $\boldsymbol{u}_i$
the outward normal vector to the facet (i.e., $n-1$-dimensional face) \sloppy
opposite to $\boldsymbol{v}_i$,
and by $\theta_{i \, j}$ the dihedral angle between two facets
with normal vectors $\boldsymbol{u}_i$ and $\boldsymbol{u}_j$
(i.e., two facets opposite to $\boldsymbol{v}_i$ and $\boldsymbol{v}_j$).
We also define $\theta_{i \, i} := \pi$.
Then the following relation holds:
\begin{equation} \label{eq_ipangle}
\left \langle \boldsymbol{u}_i , \boldsymbol{u}_j \right \rangle 
= - \cos \theta_{i \, j} \, .
\end{equation}

\subsection*{From dihedral angles to edge lengths}

Let $G_A := 
\left( \left \langle \boldsymbol{u}_i , \boldsymbol{u}_j \right \rangle \right)
_{i, j = 1}^{n+1}$ be an $ \left( n+1 \right) \times \left( n+1 \right) $ 
real matrix of signature $(n,1)$.
This matrix is called the {\em Gram matrix in terms of dihedral angles},
and determined independent of the position in $\mathbb{H}^n$ via \eqref{eq_ipangle}.
If $n=3$ and dihedral angles are named as in Figure~\ref{f_tetra},
then you can easily check that $G_A$ becomes as in \eqref{eq_ga}.

Let $M$ be the diagonal matrix with diagonal entries 
$\left \langle \boldsymbol{u}_i , \boldsymbol{v}_i \right \rangle$,
and let $\left( {c_A}_{i \, j}  \right)_{i, j = 1}^{n+1} := \left( \det G_A \right) G_A^{-1}$
be the cofactor matrix of $G_A$,
namely ${c_A}_{i \, j} := (-1)^{i+j} \det G_{i \, j}$, 
where $G_{i \, j}$ is a submatrix of order $n$ obtained from $G_A$ 
by removing the $i$th row and $j$th column.
Then we have $G_l = M {G_A}^{-1} M$.
So, using \eqref{eq_iplength} and 
$\left \langle \boldsymbol{u}_i , \boldsymbol{v}_i \right \rangle = 
- \sqrt{ ( \det G_A) / {c_A}_{i \, i} } $ for any $i$, we have
\begin{align}
\cosh \rho_{i \, j} 
	&= -\sqrt{ \frac{\det G_A}{ {c_A}_{i \, i} }}  \, \frac{{c_A}_{i \, j}}{\det G_A} \, 
		\sqrt{ \frac{\det G_A}{ {c_A}_{j \, j} }} \notag\\
	&= \frac{{c_A}_{i \, j}}{ \sqrt{ {c_A}_{i \, i} {c_A}_{j \, j} } } \, . \label{eq_a2l}
\end{align}

\subsection*{From edge lengths to dihedral angles}

The method is just the same as in the previous subsection.
Let $G_l := 
\left( \left \langle \boldsymbol{v}_i , \boldsymbol{v}_j \right \rangle \right)
_{i, j = 1}^{n+1}$ be an $ \left( n+1 \right) \times \left( n+1 \right) $
real matrix of signature $(n,1)$.
This matrix is called the {\em Gram matrix in terms of edge lengths},
and determined independent of the position in $\mathbb{H}^n$ via \eqref{eq_iplength}.
If $n=3$ and edge lengths are named as in Figure~\ref{f_tetra},
then you can easily check that $G_l$ becomes as in \eqref{eq_gl}.

Let $\left( {c_l}_{i \, j}  \right)_{i, j = 1}^{n+1} := \left( \det G_l \right) G_l^{-1}$
be the cofactor matrix of $G_l$.
Then, using \eqref{eq_ipangle} and 
$\left \langle \boldsymbol{u}_i , \boldsymbol{v}_i \right \rangle = 
- \sqrt{ ( \det G_l) / {c_l}_{i \, i} } $ for any $i$ in this case,
we have the following relationship between $\theta_{i \, j}$ and edge lengths:
\begin{equation}
\cos \theta_{i \, j} = \frac{{c_l}_{i \, j}}{ \sqrt{ {c_l}_{i \, i} {c_l}_{j \, j} } } \, . \label{eq_l2a} 
\end{equation}
%
%
%
%
%
%Remark 3.1
%
%
%
\begin{rem}
We already know a volume formula Theorem~\ref{thm_va} in terms of dihedral angles.
So combine it with \eqref{eq_l2a}
and we can theoretically obtain a formula in terms of edge lengths.
However the formula in Theorem~\ref{thm_vl} does not come from this procedure,
and its proof given in the next section is also independent of Theorem~\ref{thm_va}.
\end{rem}

%
%
%
%
%Section 4
%
%
%
%
\section{Proof of the main theorem} \label{s_proof}

We start this section with checking the details of $z_-$ and $z_+$,
most of which are already known in \cite{my}.
%
%
%
%
%Proposition 4.1
%
%
%
%
\vspace{\medskipamount} \begin{prop} \label{prop_z}
For each\/ $j = ``+" \text{ or\/ } ``-"$,
we have\/
$\displaystyle{ z_j \left. \frac{\partial U}{\partial z} \right \rvert_{z=z_j} \equiv 0 }$
modulo $2 \, \pi \, \sqrt{-1}$.
\end{prop}
\begin{proof}
Using $\displaystyle{ \frac{ \mathrm{d}}{\mathrm{d} z} \li (z) = - \frac{1}{z} \, \log(1 - z) }$,
we have
\begin{align*}
z \, \frac{\partial U}{\partial z}
	&= -\log (1 - z) - 
		\log (1 - z\,{a_1}\,{a_2}\,{a_4}\,{a_5}) - 
		\log (1 - z\,{a_1}\,{a_3}\,{a_4}\,{a_6})\\
	&\qquad - \log (1 - z\,{a_2}\,{a_3}\,{a_5}\,{a_6}) + 
		\log (1 + z\,{a_1}\,{a_2}\,{a_3}) + 
		\log (1 + z\,{a_1}\,{a_5}\,{a_6})\\
	&\qquad + \log (1 + z\,{a_2}\,{a_4}\,{a_6}) + 
		\log (1 + z\,{a_3}\,{a_4}\,{a_5}) \, .
\end{align*}
So the following are equivalent:
\begin{align*}
\lefteqn{ \Re ( z \, \frac{\partial U}{\partial z} ) = 0 }\\
	&\quad \Longleftrightarrow 
		\log \left \lvert \frac{
		(1 + z\,{a_1}\,{a_2}\,{a_3}) \, (1 + z\,{a_1}\,{a_5}\,{a_6}) \, 
		(1 + z\,{a_2}\,{a_4}\,{a_6}) \, (1 + z\,{a_3}\,{a_4}\,{a_5})}{
		(1 - z) \, (1 - z\,{a_1}\,{a_2}\,{a_4}\,{a_5}) \, 
		(1 - z\,{a_1}\,{a_3}\,{a_4}\,{a_6}) \, (1 - z\,{a_2}\,{a_3}\,{a_5}\,{a_6})} 
		\right \rvert = 0 \\
	&\quad \Longleftrightarrow 
		(1 + z\,{a_1}\,{a_2}\,{a_3}) \, (1 + z\,{a_1}\,{a_5}\,{a_6}) \, 
		(1 + z\,{a_2}\,{a_4}\,{a_6}) \, (1 + z\,{a_3}\,{a_4}\,{a_5})\\
	&\qquad \qquad \pm (1 - z) \, (1 - z\,{a_1}\,{a_2}\,{a_4}\,{a_5}) \, 
		(1 - z\,{a_1}\,{a_3}\,{a_4}\,{a_6}) \, (1 - z\,{a_2}\,{a_3}\,{a_5}\,{a_6}) = 0 \, .
\end{align*}
When we choose ``$-$" part of the above equation, the following equation holds:
\begin{align*}
\lefteqn{ (1 + z\,{a_1}\,{a_2}\,{a_3}) \, (1 + z\,{a_1}\,{a_5}\,{a_6}) \, 
		(1 + z\,{a_2}\,{a_4}\,{a_6}) \, (1 + z\,{a_3}\,{a_4}\,{a_5}) }\\
	&\quad - (1 - z) \, (1 - z\,{a_1}\,{a_2}\,{a_4}\,{a_5}) \, 
		(1 - z\,{a_1}\,{a_3}\,{a_4}\,{a_6}) \, (1 - z\,{a_2}\,{a_3}\,{a_5}\,{a_6})\\
	&\qquad = z \left( q_2 \, z^2 + q_1 \, z + q_0 \right) .
\end{align*}
Thus $z_j$ is a solution of $\displaystyle{ \Re ( z \, \frac{\partial U}{\partial z} )= 0 }$,
namely $\displaystyle{ \Re ( z_j \left. \frac{\partial U}{\partial z} \right \rvert_{z=z_j} ) = 0 }$.
\vspace{\smallskipamount}

On the other hand, by the definition of $\log$, 
the following holds modulo $2 \, \pi$:
\begin{align*}
\lefteqn{ \Im ( z_j \left. \frac{\partial U}{\partial z} \right \rvert_{z=z_j} ) }\\
	&\quad \equiv \arg \frac{
		(1 + z_j\,{a_1}\,{a_2}\,{a_3}) \, (1 + z_j\,{a_1}\,{a_5}\,{a_6}) \, 
		(1 + z_j\,{a_2}\,{a_4}\,{a_6}) \, (1 + z_j\,{a_3}\,{a_4}\,{a_5})}{
		(1 - z_j) \, (1 - z_j\,{a_1}\,{a_2}\,{a_4}\,{a_5}) \, 
		(1 - z_j\,{a_1}\,{a_3}\,{a_4}\,{a_6}) \, (1 - z_j\,{a_2}\,{a_3}\,{a_5}\,{a_6})}\\
	&\quad = 0 \, ,
\end{align*}
where $\Im$ means the imaginary part.
We have thus proved Proposition~\ref{prop_z}.
\end{proof}

Proposition~\ref{prop_z} means that 
$\displaystyle{ z_j \left. \frac{\partial U}{\partial z} \right \rvert_{z=z_j} }$
can take discrete values.
But, by the definition,  
$\displaystyle{ z_j \left. \frac{\partial U}{\partial z} \right \rvert_{z=z_j} }$
is an analytic function in terms of $a_1, a_2, \ldots, a_6$.
So $\displaystyle{ \frac{\partial}{\partial a_i} 
	\left( z_j \left. \frac{\partial U}{\partial z} \right \rvert_{z=z_j} \right) }$ must be zero.
Thus we have the following corollary:
%
%
%
%
%Corollary 4.2
%
%
%
%
\begin{cor}\label{cor_der}
For any\/ $i = 1, 2, \ldots , 6$, we have
$\displaystyle{ \frac{\partial}{\partial a_i} 
\left( z_j \left. \frac{\partial U}{\partial z} \right \rvert_{z=z_j} \right) = 0 }$. \qed
\end{cor}

Next we prepare the following proposition:
%
%
%
%
%Proposition 4.3
%
%
%
%
\begin{prop} \label{prop_dv}
For each\/ $i = 1, 2, \ldots , 6$, we have\/
$\displaystyle{ 2 \, \frac{\partial V_l}{\partial l_i} \equiv A_i  }$
modulo $\pi$.
\end{prop}
\begin{proof}
Here we only show that $\displaystyle{ 2 \, \frac{\partial V_l}{\partial l_1} \equiv A_1}$ 
modulo $\pi$.
Other cases can be shown similarly.
Using Corollary~\ref{cor_der},
the following holds for each $j = ``+" \text{ or } ``-"$:
\begin{align*}
\lefteqn{\frac{\partial}{\partial l_1} \left( U (a_1, a_2, a_3, a_4, a_5, a_6, z_j)  
	- z_j \left. \frac{\partial U}{\partial z} \right \rvert_{z = z_j} \, \log z_j \right) }\\
		&\quad = \frac{\partial}{\partial a_4} \left( 
			U (a_1, a_2, a_3, a_4, a_5, a_6, z_j)  
			- z_j \left. \frac{\partial U}{\partial z} \right \rvert_{z = z_j} \, \log z_j  \right) 
			\, \frac{\mathrm{d} a_4}{\mathrm{d} l_1}\\
		&\quad = a_4 \, \Biggl \{ 
			\left. \frac{\partial U}{\partial a_4} \right \rvert_{z = z_j} + 
			\left. \frac{\partial U}{\partial z} \right \rvert_{z = z_j} \, 
			\frac{\partial z_j}{\partial a_4}\\
		&\quad \qquad -
			\left( \log z_j  \right) \frac{\partial}{\partial a_4} \left( 
			z_j \left. \frac{\partial U}{\partial z} \right \rvert_{z = z_j} \right)
			-  \left( z_j \left. \frac{\partial U}{\partial z} \right \rvert_{z = z_j} \right)
			\left( \frac{1}{z_j} \, \frac{\partial z_j}{\partial a_4} \right) \Biggr \}\\
		&\quad = a_4 \left. \frac{\partial U}{\partial a_4} \right \rvert_{z = z_j}\\
		&\quad \equiv \log \frac{\varphi (z_j)}{\psi (z_j)} \pmod {2 \, \pi \, \sqrt{-1} } \, ,
\end{align*}
where $\varphi (z)$ and $\psi (z)$ are defined as follows:
\begin{align*}
\varphi (z) &:= \left( 1 + z\,{a_2}\,{a_4}\,{a_6} \right) \, 
	\left( 1 + z\,{a_3}\,{a_4}\,{a_5} \right) ,\\
\psi (z) &:= \left( 1 - z\,{a_1}\,{a_2}\,{a_4}\,{a_5} \right) \,
	\left( 1 - z\,{a_1}\,{a_3}\,{a_4}\,{a_6} \right) \, .
\end{align*}
Thus the following holds:
\begin{align*}
2 \, \frac{\partial V_l}{\partial l_1}
	&= \frac{\sqrt{-1}}{2} \, \frac{\partial}{\partial l_1} \left\{
		\left( U (a_1, a_2, a_3, a_4, a_5, a_6, z_-) - 
		z_- \left. \frac{\partial U}{\partial z} \right \rvert_{z = z_-} \, \log z_- \right) \right.\\
	&\qquad \qquad \qquad \left. - \left( U (a_1, a_2, a_3, a_4, a_5, a_6, z_+) - 
		z_+	\left. \frac{\partial U}{\partial z} \right \rvert_{z = z_+} \, \log z_+ \right) \right\}\\
	&\equiv \frac{\sqrt{-1}}{2} \left( \log \frac{\varphi (z_-)}{\psi (z_-)} 
		- \log \frac{\varphi (z_+)}{\psi (z_+)} \right)  \pmod \pi \\
	&\equiv \frac{\sqrt{-1}}{2} \log \frac{\varphi (z_-) \, \psi (z_+)}
		{\varphi (z_+) \, \psi (z_-)} \pmod {2 \, \pi} \, .
\end{align*}
Now 
\begin{align*}
\varphi (z_-) \, \psi (z_+) 
	&= \frac{4 \, {a_1}\,{a_2}^2\,{a_3}^2\,{{a_4}}^4\,{a_5}^2\,{a_6}^2}{{q_2}^2}
		\left( {a_1}\,{a_2} + {a_3} \right) \left( {a_2} + {a_1}\,{a_3} \right) \times \\
	&\qquad \left( {a_1}\,{a_5} + {a_6} \right) \left( {a_5} + {a_1}\,{a_6} \right) 
		\left( {c_l}_{1\,2} - \sqrt{-1} 
		\sqrt{ {c_l}_{1\,1} \, {c_l}_{2\,2} - {c_l}_{1\,2}^2 } \right) ,\\
\varphi (z_+) \, \psi (z_-) 
	&= \overline{\varphi (z_-) \, \psi (z_+)} \, ,
\end{align*}
so we have
\begin{align*}
2 \, \frac{\partial V_l}{\partial l_1}
	&\equiv \frac{\sqrt{-1}}{2} \log \frac{\varphi (z_-) \, \psi (z_+)}
		{\, \overline{ \varphi (z_-) \, \psi (z_+) } \,}  \pmod \pi \\
	&= - \arg \varphi (z_-) \, \psi (z_+)\\
	&= \arctan \frac{\sqrt{ {c_l}_{1\,1} \, {c_l}_{2\,2} - {c_l}_{1\,2}^2 }}{{c_l}_{1\,2}}\\
	&= A_1 \, ,
\end{align*}
where we use \eqref{eq_l2a} at the last equality.
Thus we have proved  Proposition~\ref{prop_dv}.
\end{proof}

We finally recall the so-called Schl\"afli's differential formula:
%
%
%
%
%Proposition Schl\"afli
%
%
%
%
\begin{prop}[see, for example, \cite{mi}]
Let\/ $T$ be a hyperbolic tetrahedron with dihedral angle\/ $A_i$ and edge length\/ $l_i$.
Then the following equation holds:
\begin{equation*}
\mathrm{d} \vol (T) = - \frac{1}{2} \, \sum_{i=1}^6 l_i \, \mathrm{d} A_i \, .
\end{equation*}
\end{prop}

\begin{proof}[Proof of Theorem~\ref{thm_vl}]
Let $W = W (T)$ be the right-hand side of the formula.
Since $W (T)$ is represented by the edge lengths symmetrically,
and the edge lengths can be represented by the dihedral angles symmetrically as in \eqref{eq_a2l},
all we need to show are, by the Schl\"afli's differential formula, the following two facts:
\begin{enumerate}
\item \label{i_dif} $
	\displaystyle{\frac{\partial W}{\partial A_1} = - \frac{l_1}{2} \, ,}$
\item \label{i_val} 
	At some value $T$ of a hyperbolic tetrahedron, $W (T)$ is actually its volume. 
\end{enumerate}

\vspace{\bigskipamount}

Using Proposition~\ref{prop_dv}, the former one can be shown as follows:

\begin{align*}
\frac{\partial W}{\partial A_1} 
&= \frac{\partial V_l}{\partial A_1}  - 
	\sum_{i=1}^6 \frac{\partial}{\partial A_1} \left( l_i \frac{\partial V_l}{\partial l_i} \right)\\
&= \frac{\partial V_l}{\partial A_1}  - 
	\sum_{i=1}^6 \left( \frac{\partial l_i}{\partial A_1} \frac{\partial V_l}{\partial l_i} 
	+ l_i \frac{\partial}{\partial A_1} \frac{\partial V_l}{\partial l_i} \right)\\
&= \frac{\partial V_l}{\partial A_1}  
	- \sum_{i=1}^6 \frac{\partial V_l}{\partial l_i} \frac{\partial l_i}{\partial A_1} 
	- \sum_{i=1}^6  \frac{l_i}{2} \frac{\partial}{\partial A_1} \left( A_i + k_i \, \pi \right) 
	\quad ( k_i \in \mathbb{Z} )\\
&= - \frac{l_1}{2} \, .
\end{align*}
Since both $W (T)$ and $\vol (T)$ are analytic function with variables $l_1 , l_2 , \ldots , l_6$,
the result obtained above means that the difference between $W (T)$ and $\vol (T)$ is a constant
independent of $l_1 , l_2 , \ldots , l_6$.

\vspace{\medskipamount}

For the proof of the later one, 
let $T_\rho$ be a hyperbolic regular tetrahedron with all edge lengths being $\rho$,
and we consider a sequence of such tetrahedra shrinking to a point.
Then  clearly $\displaystyle{\lim_{\rho \rightarrow 0} \vol (T_\rho) = 0}$.

On the other hand, we have the following equations in this case:
\begin{align*}
z_\pm 
	&= \frac{3 \left( \cosh \rho + 1 \right) \pm 
		\sqrt{-1} \sqrt{\left( \cosh \rho - 1 \right) \left( 3 \, \cosh \rho + 1 \right)}}
		{2 \, e^{4 \, \rho} \left( 2 \, \cosh \rho - \sinh \rho + 1 \right)} \, ,\\
U 
	&= \li (z_j) + 3 \, \li (e^{4 \, \rho} z_j) - 4 \, \li (e^{3 \, \rho} z_j) \, ,\\
z_j \left. \frac{\partial U}{\partial z} \right \rvert_{z = z_j} 
	&= -\sqrt{-1} \left\{ \arg (1-z_j) + 3 \, \arg (1- e^{4 \, \rho} z_j) - 
		4 \, \arg (1- e^{3 \, \rho} z_j ) \right\} ,\\
\frac{\partial V_l}{\partial \rho} 
	&= \frac{1}{2} \left( \arctan 
		\frac{\sqrt{\left( \cosh \rho + 1 \right) \left( 3 \, \cosh \rho + 1 \right)}}{\cosh \rho} 
		+ k \, \pi  \right) \quad ( k \in \mathbb{Z}) \, ,
\end{align*}
where the third equation follows from Proposition~\ref{prop_z},
and the last one from Proposition~\ref{prop_dv}.
Thus the following equation holds:

\begin{align*}
W (T_\rho) 
&= \frac{\sqrt{-1}}{4} \Bigl\{ 
	\left\{ \li (z_-) - \li (z_+) \right\}
	+ 3 \left\{ \li (e^{4 \, \rho} z_-) - \li (e^{4 \, \rho} z_+) \right\}\\
&\quad - 4 \left\{ \li (e^{3 \, \rho} z_-) - \li (e^{3 \, \rho} z_+) \right\}\\
&\quad + \sqrt{-1} \left\{ \arg (1-z_-) + 3 \arg (1- e^{4 \, \rho} z_-) - 
	4 \arg (1- e^{3 \, \rho} z_- ) \right\} \log z_- \\
&\quad - \sqrt{-1} \left\{ \arg (1-z_+ ) + 3 \arg (1- e^{4 \, \rho} z_+ ) - 
	4 \arg (1- e^{3 \, \rho} z_+ ) \right\} \log z_+ \Bigr\}\\
&\quad - 3 \, \rho \left(
	\arctan \frac{\sqrt{\left( \cosh \rho + 1 \right) \left( 3 \, \cosh \rho + 1 \right)}}
		{\cosh \rho}  + k \, \pi \right)\\
&= - \frac{1}{2} \Bigl\{ 
	\Im \left\{ \li (z_-) + 3 \li (e^{4 \, \rho} z_-)  - 4 \li (e^{3 \, \rho} z_-) \right\}\\
&\qquad + \left\{ \arg (1-z_-) + 3 \arg (1- e^{4 \, \rho} z_-) - 
	4 \arg (1- e^{3 \, \rho} z_- ) \right\} \log \lvert z_- \rvert \Bigr\}\\
&\qquad - 3 \, \rho \left(
	\arctan \frac{\sqrt{\left( \cosh \rho + 1 \right) \left( 3 \, \cosh \rho + 1 \right)}}
		{\cosh \rho} + k \, \pi \right) ,
\end{align*}
where we used the fact $z_+ = \overline{z_-}$ in the last equation.

When $\rho$ converges to $0$, then $z_-$ converges to $1$, a branch point of $\li (z)$.
However this convergence is not spiral 
because $\Im z_- \ne 0$ for any  $\rho > 0$. So we have
\begin{equation*}
\lim_{\rho \rightarrow 0} \Im 
\left\{ \li (z_-) + 3 \li (e^{4 \, \rho} z_-)  - 4 \li (e^{3 \, \rho} z_-) \right\} = 0 \, .
\end{equation*}
Furthermore we have 
\begin{align*}
\arg (1-z_-) 
	&= \arctan
		\frac{{-\sqrt{3\,e^{2\,\rho } + 2\,e^{\rho } + 3}}}
		{2\,e^{5\,\rho } + 6\,e^{4\,\rho } + 12\,e^{3\,\rho } + 
		12\,e^{2\,\rho } + 9\,e^{\rho } + 3 } \, ,\\
\arg (1- e^{4 \, \rho} z_-) 
	&= \arctan 
		\frac{{\sqrt{3\,e^{2\,\rho } + 2\,e^{\rho } + 3}}}{e^{\rho } + 3} \, ,\\
\arg (1- e^{3 \, \rho} z_-) 
	&= \arctan 
		\frac{{ -\sqrt{3\,e^{2\,\rho } + 2\,e^{\rho } + 3}}}
		{2\,e^{2\,\rho } + 3\,e^{\rho } + 3} \, .
\end{align*}
Thus we conclude as follows:
\begin{equation*}
\lim_{\rho \rightarrow 0} \left\{ \arg (1-z_-) + 3 \arg (1- e^{4 \, \rho} z_-) - 
	4 \arg (1- e^{3 \, \rho} z_- ) \right\} \log \lvert z_- \rvert = 0 \, .
\end{equation*}
Finally we have 
\begin{equation*}
\lim_{\rho \rightarrow 0} \frac{\sqrt{\left( \cosh \rho + 1 \right) \left( 3 \, \cosh \rho + 1 \right)}}
		{\cosh \rho} = 2 \sqrt{2} 
\end{equation*}
and $\arctan 2 \sqrt{2} $ is actually the dihedral angle of the Euclidean regular tetrahedron.

Thus we have $\displaystyle{ \lim_{\rho \rightarrow 0} W (\rho) = 0 }$ and 
finished the proof of Theorem~\ref{thm_vl}.
\end{proof}

%%%%%%%%%%

\end{document}